# Complete moment and integral convergence for sums of negatively associated random variables

**Han-Ying Liang**[*]

*Department of Mathematics, Tongji University*
*Shanghai 200092, P. R. China*
*e-mail:* hyliang83@yahoo.com

**Deli Li**[†]

*Department of Mathematical Sciences, Lakehead University*
*Thunder Bay, Ontario, Canada P7B 5E1*
*e-mail:* dli@lakeheadu.ca

**Andrew Rosalsky**

*Department of Statistics, University of Florida*
*Gainesville, FL 32611, USA*
*e-mail:* rosalsky@stat.ufl.edu

**Abstract:** For a sequence of identically distributed negatively associated random variables $\{X_n; \; n \geq 1\}$ with partial sums $S_n = \sum_{i=1}^{n} X_i$, $n \geq 1$, refinements are presented of the classical Baum-Katz and Lai complete convergence theorems. More specifically, necessary and sufficient moment conditions are provided for complete moment convergence of the form

$$\sum_{n \geq n_0} n^{r-2-\frac{1}{pq}} a_n E \left( \max_{1 \leq k \leq n} |S_k|^{\frac{1}{q}} - \epsilon b_n^{\frac{1}{pq}} \right)^+ < \infty$$

to hold where $r > 1$, $q > 0$ and either $n_0 = 1$, $0 < p < 2$, $a_n = 1$, $b_n = n$ or $n_0 = 3$, $p = 2$, $a_n = (\log n)^{-\frac{1}{2q}}$, $b_n = n \log n$. These results extend results of Chow (1988) and Li and Spătaru (2005) from the independent and identically distributed case to the identically distributed negatively associated setting. The complete moment convergence is also shown to be equivalent to a form of complete integral convergence.

**AMS 2000 subject classifications:** Primary 60G50; secondary 60F15.
**Keywords and phrases:** Baum-Katz's law, Lai's law, complete moment convergence, complete integral convergence, convergence rate of tail probabilities, sums of identically distributed and negatively associated random variables.

[*]The research of Han-Ying Liang was supported by the National Natural Science Foundation of China (10571136).

[†]The research of Deli Li was supported by a grant from the Natural Sciences and Engineering Research Council of Canada.

## 1. Introduction

Let $\{X_n; n \geq 1\}$ be a sequence of random variables and, as usual, set $S_n = X_1 + X_2 + \cdots + X_n$, $n \geq 1$. When $\{X, X_n; n \geq 1\}$ are independent and identically distributed (i.i.d.), Baum and Katz (1965) proved the following remarkable result concerning the convergence rate of the tail probabilities $P\left(|S_n| > \epsilon n^{1/p}\right)$ for arbitrary $\epsilon > 0$.

**Theorem 1.1** (Baum and Katz (1965))**.** *Let* $0 < p < 2$ *and* $\gamma \geq p$ *. Then*

$$\sum_{n \geq 1} n^{\frac{\gamma}{p} - 2} P\left(|S_n| > \epsilon n^{\frac{1}{p}}\right) < \infty \quad \text{for all} \quad \epsilon > 0$$

*if and only if* $E|X|^\gamma < \infty$, *where* $EX = 0$ *whenever* $1 \leq p < 2$.

This is what is known as the Baum-Katz law. Some special cases of the Baum-Katz law are important well-known theorems. The special case $p = 1$ and $\gamma = 2$ is the famous Hsu-Robbins-Erdös complete convergence theorem (see Hsu and Robbins (1947) for the sufficiency half and Erdös (1949, 1950) for the necessity half) whereas the special case $p = \gamma = 1$ is the celebrated theorem of Spitzer (1956). Moreover, there is an interesting and substantial literature of investigation apropos of extending the Hsu-Robbins-Erdös, Spitzer, and Baum-Katz theorems along a variety of different paths. One of these extensions is due to Chow (1988) who established the following refinement which is a complete moment convergence result for sums of i.i.d. random variables.

**Theorem 1.2** (Chow (1988))**.** *Let* $EX = 0$, *and let* $1 \leq p < 2$ *and* $\gamma \geq p$. *Suppose that* $E[|X|^\gamma + |X| \log(1 + |X|)] < \infty$. *Then*

$$\sum_{n=1}^{\infty} n^{\frac{\gamma}{p} - 2 - \frac{1}{p}} E\left(|S_n| - \epsilon n^{\frac{1}{p}}\right)^+ < \infty \quad \text{for all} \quad \epsilon > 0.$$

By applying the Hoffmann-Jørgensen (1974) inequality, Li and Spătaru (2005) obtained the following generalization of Theorem 1.2, also for sums of i.i.d. random variables. (It will be transparent from Remark 1.1 below that Theorem 1.3 does generalize Theorem 1.2.) Theorem 1.3 is a refinement of Theorem 1.1.

**Theorem 1.3** (Li and Spătaru (2005), Theorem 2)**.** *Let* $EX = 0$, *and let* $0 < p < 2$, $r \geq 1$ *and* $q > 0$. *Set*

$$f(x) = \sum_{n \geq 1} n^{r-2} P\left(|S_n| > x n^{\frac{1}{p}}\right), \quad x > 0.$$

*Then the following are equivalent:*

(i) $\displaystyle\int_\epsilon^\infty f(x^q)\, dx < \infty$ *for all* $\epsilon > 0$;

(ii) $\begin{cases} E|X|^{\frac{1}{q}} < \infty & \text{if } q < \frac{1}{pr}, \\ E[|X|^{pr} \log^+ |X|] < \infty & \text{if } q = \frac{1}{pr}, \\ E|X|^{pr} < \infty & \text{if } q > \frac{1}{pr}. \end{cases}$

**Remark 1.1.** *It is natural to refer to Theorem 1.3 as being a complete integral convergence theorem. We will now verify that for all $\epsilon > 0$,*

$$\int_\epsilon^\infty f(x^q)\,dx = \sum_{n \geq 1} n^{r-2-\frac{1}{pq}} E\left(|S_n|^{\frac{1}{q}} - \epsilon n^{\frac{1}{pq}}\right)^+. \tag{1.1}$$

*Consequently, Theorem 1.3 is also a complete moment convergence theorem and when $q = 1$, $r > 1$, $1 \leq p < 2$, the implication ((ii) $\Rightarrow$ (i)) in Theorem 1.3 reduces to Theorem 1.2 by taking $\gamma = pr$. To verify (1.1), observe for $\epsilon > 0$ that*

$$\begin{aligned}
\int_\epsilon^\infty f(x^q)\,dx &= \sum_{n \geq 1} n^{r-2} \int_\epsilon^\infty P\left(|S_n| > x^q n^{\frac{1}{p}}\right) dx \\
&= \sum_{n \geq 1} n^{r-2-\frac{1}{pq}} \left(\frac{1}{q} \int_{\epsilon^q n^{1/p}}^\infty t^{\frac{1}{q}-1} P(|S_n| > t)\,dt\right) \\
&\quad \left(by\ letting\ t = x^q n^{\frac{1}{p}}\right)
\end{aligned} \tag{1.2}$$

*and that*

$$\begin{aligned}
&E\left(|S_n|^{\frac{1}{q}} I\left(|S_n| > \epsilon^q n^{1/p}\right)\right) \\
&= \frac{1}{q} \int_0^\infty t^{\frac{1}{q}-1} P\left(|S_n| I\left(|S_n| > \epsilon^q n^{\frac{1}{p}}\right) > t\right) dt \\
&= \frac{1}{q} \int_0^{\epsilon^q n^{1/p}} t^{\frac{1}{q}-1} P\left(|S_n| I\left(|S_n| > \epsilon^q n^{\frac{1}{p}}\right) > t\right) dt \\
&\quad + \frac{1}{q} \int_{\epsilon^q n^{1/p}}^\infty t^{\frac{1}{q}-1} P\left(|S_n| I\left(|S_n| > \epsilon^q n^{\frac{1}{p}}\right) > t\right) dt \\
&= \frac{1}{q} \int_0^{\epsilon^q n^{1/p}} t^{\frac{1}{q}-1} P\left(|S_n| > \epsilon^q n^{\frac{1}{p}}\right) dt + \frac{1}{q} \int_{\epsilon^q n^{1/p}}^\infty t^{\frac{1}{q}-1} P(|S_n| > t)\,dt \\
&= \epsilon n^{\frac{1}{pq}} P\left(|S_n| > \epsilon^q n^{\frac{1}{p}}\right) + \frac{1}{q} \int_{\epsilon^q n^{1/p}}^\infty t^{\frac{1}{q}-1} P(|S_n| > t)\,dt.
\end{aligned}$$

*Thus,*

$$\begin{aligned}
&\frac{1}{q} \int_{\epsilon^q n^{1/p}}^\infty t^{\frac{1}{q}-1} P(|S_n| > t)\,dt \\
&= E\left(|S_n|^{\frac{1}{q}} I\left(|S_n| > \epsilon^q n^{\frac{1}{p}}\right)\right) - \epsilon n^{\frac{1}{pq}} P\left(|S_n| > \epsilon^q n^{\frac{1}{p}}\right) \\
&= E\left(\left(|S_n|^{\frac{1}{q}} - \epsilon n^{\frac{1}{pq}}\right) I\left(|S_n| > \epsilon^q n^{\frac{1}{p}}\right)\right) \\
&= E\left(|S_n|^{\frac{1}{q}} - \epsilon n^{\frac{1}{pq}}\right)^+
\end{aligned}$$

*which, when combined with (1.2), yields (1.1).*

　Additionally, and again for the i.i.d. case, Li and Spătaru (2005) refined Lai's (1974) law by the following complete integral convergence theorem.

**Theorem 1.4** (Li and Spătaru (2005), Theorem 4). *Let $r > 1$ and $q > 0$. Set*

$$f(x) = \sum_{n \geq 2} n^{r-2} P(|S_n| > x\sqrt{n \log n}), \quad x > 0.$$

(i) *If $\int_\epsilon^\infty f(x^q)dx < \infty$ for some $\epsilon > 0$, then $EX = 0$ and*

$$\begin{cases} E|X|^{\frac{1}{q}} < \infty & \text{if } q \leq \frac{1}{2r}, \\ E[|X|^{2r}(\log^+ |X|)^{-r}] < \infty & \text{if } q > \frac{1}{2r}. \end{cases} \quad (1.3)$$

(ii) *If (1.3) holds, $EX^2 = \sigma^2$, and $EX = 0$, then $\int_\epsilon^\infty f(x^q)dx < \infty$, $\epsilon > \left(\sigma\sqrt{2r-2}\right)^{\frac{1}{q}}$.*

**Remark 1.2.** *Apropos of the function $f(x)$ in Theorem 1.4, for all $\epsilon > 0$*

$$\int_\epsilon^\infty f(x^q)\,dx = \sum_{n \geq 2} n^{r-2-\frac{1}{2q}}(\log n)^{-\frac{1}{2q}} E\left(|S_n|^{\frac{1}{q}} - \epsilon(n \log n)^{\frac{1}{2q}}\right)^+$$

*by arguing as in Remark 1.1, mutatis mutandis. Consequently, Theorem 1.4 is also a complete moment convergence theorem.*

In the current work, we establish in Theorems 2.1 and 2.2 complete moment and integral convergence theorems which are new versions of Theorems 1.3 and 1.4, respectively where the previous i.i.d. hypothesis is replaced by the assumption that $\{X, X_n; n \geq 1\}$ are identically distributed and negatively associated and where $S_n$ is replaced by $\max_{1 \leq k \leq n} |S_k|$. Theorem 2.1 (resp., Theorem 2.2) thus concerns the convergence rate of the tail probabilities $P\left(\max_{1 \leq k \leq n} |S_k| > xn^{1/p}\right)$ (resp., $P\left(\max_{1 \leq k \leq n} |S_k| > x\sqrt{n \log n}\right)$).

We now recall the definition of negatively associated random variables.

**Definition:** *A finite family of random variables $\{X_i; 1 \leq i \leq n\}$ is said to be negatively associated (NA) if, for every pair of disjoint subsets $A$ and $B$ of $\{1, 2, \cdots, n\}$, we have*

$$\text{Cov}(f_1(X_i, i \in A), f_2(X_j, j \in B)) \leq 0$$

*whenever $f_1$ and $f_2$ are coordinatewise increasing and the covariance exists. An infinite family of random variables is NA if every finite subfamily is NA.*

Alam and Saxena (1981) introduced the notion of negative association. Many well-known families of multivariate distributions have the NA property. A prominent example is the multivariate normal distribution with negatively correlated components. See Joag-Dev and Proschan (1983) for a listing and verification. In many stochastic models, an independence assumption among the random variables in the model is not a reasonable assumption since the random variables may be "repelling" in the sense that increases in any of the random variables often correspond to decreases in the others. Thus, the NA assumption is often more suitable than the classical assumption of independence and the concept of negative association is useful for a diversity of applications in multivariate analysis. Basic properties of NA families are developed in Joag-Dev and Proschan (1983).

The main classes of limit theorems and moment inequalities for partial sums of independent random variables have counterparts for NA random variables. Some references for these NA counterparts are listed as follows:

- Central Limit Theorem (Newman (1984), Roussas (1994), Su and Chi (1998))
- Three-Series Criterion (Matula (1992))
- Weak Law of Large Numbers (Chi and Su (1997))
- Moment Inequalities (Su et al. (1997), Liu et al. (1999), Shao (2000))
- Law of the Iterated Logarithm (Shao and Su (1999))
- Complete Convergence Theorem (Liang and Su (1999), Liang (2000), Baek et al. (2003))
- Strong Law of Large Numbers (Liang and Baek (2006)).

The concept of complete moment convergence was introduced by Chow (1988) for the case of i.i.d. summands. Extensions for sums of Banach space valued random elements were investigated by Wang and Su (2004), Wang et al. (2005), Chen (2006), Guo and Xu (2006), Rosalsky et al. (2006), and Ye and Zhu (2007). Wang and Zhao (2006) studied complete moment convergence for sums of NA random variables but their results and ours do not entail each other. Precise asymptotics for complete moment convergence were studied by Jiang and Zhang (2006), Liu and Lin (2006), Ye et al. (2007), and Fu and Zhang (2008) for i.i.d. summands and by Jiang and Zhang (2005) for a stationary sequence of NA summands.

The concept of complete integral convergence was introduced by Spătaru (1990) for the case of i.i.d. summands wherein he strengthened the Hsu-Robbins-Erdös theorem. Further work on complete integral convergence with i.i.d. summands was carried out by Li and Spătaru (2005) and by Spătaru (2005, 2006).

Throughout this paper, the symbol $C$ will denote a generic positive constant whose value may change from one place to another.

The layout of the paper is as follows. Our main results regarding the refinement of the Baum-Katz and Lai laws for NA random variables, Theorems 2.1 and 2.2, are presented in Section 2. Their proofs will be provided in Section 3. The method used for proving our main results is different from that of Li and Spătaru (2005).

## 2. The Main Results

In this section, let $\{X, X_n; n \geq 1\}$ be a sequence of identically distributed NA random variables, and set $S_n = \sum_{i=1}^{n} X_i, n \geq 1$. The statements of the main results on complete moment and integral convergence follow.

**Theorem 2.1.** *Let* $0 < p < 2$, $r > 1$, $q > 0$, *and let*

$$f(x) = \sum_{n \geq 1} n^{r-2} P\left(\max_{1 \leq k \leq n} |S_k| > xn^{\frac{1}{p}}\right), \quad x > 0.$$

*Then the following are equivalent:*

(i) $\sum_{n \geq 1} n^{r-2-\frac{1}{pq}} E\left(\max_{1 \leq k \leq n} |S_k|^{\frac{1}{q}} - \epsilon n^{\frac{1}{pq}}\right)^+ < \infty$ *for all* $\epsilon > 0$;

(ii) $\int_{\epsilon}^{\infty} f(x^q)\, dx < \infty$ *for all* $\epsilon > 0$;

(iii) $\begin{cases} E|X|^{\frac{1}{q}} < \infty & \text{if } q < \frac{1}{pr}, \\ E[|X|^{pr} \log^+ |X|] < \infty & \text{if } q = \frac{1}{pr}, \\ E|X|^{pr} < \infty & \text{if } q > \frac{1}{pr}, \\ EX = 0 & \text{when } 1 \leq p < 2. \end{cases}$

**Remark 2.1.** *Some comments comparing Theorem 2.1 with the main result of Wang and Zhao (2006) are in order. For a sequence of NA random variables $\{X_n; n \geq 1\}$ which are stochastically dominated by a random variable $X$, Wang and Zhao (2006) considered only the case $E|X|^p < \infty$ where $p > 1$ and showed that this moment condition gives a complete moment convergence result. Furthermore, they showed that the moment condition $E|X|^p < \infty$ $(p > 1)$ is a necessary condition for their complete moment convergence result if the NA sequence also stochastically majorizes $X$. The stochastic domination and majorization conditions are of course automatic if the sequence $\{X_n; n \geq 1\}$ consists of identically distributed random variables. On the other hand, in our Theorem 2.1, we allow for all possibilities relating $q$ with $pr$ in the moment conditions comprising condition (iii). Thus in the equivalence provided by Theorem 2.1, we allow for $X$ to have a finite absolute moment of any order. The Wang and Zhao (2006) result and Theorem 2.1 have different proofs.*

**Theorem 2.2.** *Let* $r > 1$, $q > 0$, *and let*

$$f(x) = \sum_{n \geq 3} n^{r-2} P\left(\max_{1 \leq k \leq n} |S_k| > x\sqrt{n \log n}\right), \quad x > 0.$$

*Then the following are equivalent:*

(i) $\sum_{n \geq 3} n^{r-2-\frac{1}{2q}} (\log n)^{-\frac{1}{2q}} E\left(\max_{1 \leq k \leq n} |S_k|^{\frac{1}{q}} - \epsilon(n \log n)^{\frac{1}{2q}}\right)^+ < \infty$
  *for some* $\epsilon > 0$;

(ii) $\int_{\epsilon}^{\infty} f(x^q)\, dx < \infty$ *for some* $\epsilon > 0$;

$$\text{(iii)} \begin{cases} E|X|^{\frac{1}{q}} < \infty & \text{if } q \leq \frac{1}{2r}, \\ E[|X|^{2r}(\log^+ |X|)^{-r}] < \infty & \text{if } q > \frac{1}{2r}, \\ EX = 0. \end{cases}$$

## 3. Proofs of the Main Results

The following four lemmas will be needed to prove the main results.

**Lemma 3.1.** *Let* $\{X_i; \ i \geq 1\}$ *be a sequence of NA random variables and let* $\{a_{ni}; \ 1 \leq i \leq n, n \geq 1\}$ *be an array of real numbers. Then there exists a constant* $A > 0$ *such that for all* $n \geq 1$ *and all* $\epsilon > 0$,

$$\frac{1}{2} \sum_{j=1}^{n} P(|a_{nj}X_j| > \epsilon) \leq (1+A) P\left(\max_{1 \leq j \leq n} |a_{nj}X_j| > \epsilon\right)$$
$$+ \sum_{j=1}^{n} P(|a_{nj}X_j| > \epsilon) P\left(\max_{1 \leq j \leq n} |a_{nj}X_j| > \epsilon\right).$$

*Moreover, if* $\lim_{n \to \infty} P(\max_{1 \leq j \leq n} |a_{nj}X_j| > \epsilon) = 0$ *for all* $\epsilon > 0$, *then there exists a constant* $C > 0$ *such that for all sufficiently large* $n$ *and all* $\epsilon > 0$,

$$\sum_{j=1}^{n} P(|a_{nj}X_j| > \epsilon) \leq C P\left(\max_{1 \leq j \leq n} |a_{nj}X_j| > \epsilon\right).$$

The proof of Lemma 3.1 can be found in the proof of Theorem 2.1 in Liang and Su (1999).

**Lemma 3.2** (Li and Spǎtaru (2005), Lemma 4). *Let* $\psi$ *be an increasing nonnegative function mapping* $[0, \infty)$ *onto* $[0, \infty)$, *with inverse function* $\psi^{-1}$, *and let* $\phi$ *be a positive nondecreasing function on* $[0, \infty)$. *Assume that* $\psi$ *and* $\phi$ *are regularly varying at infinity with respective indexes* $\alpha > 0$ *and* $\beta \geq 0$, *and define*

$$\Psi(x) = \begin{cases} 0 & \text{if } 0 \leq x \leq 1, \\ \int_1^x \phi(y)/\psi(y) dy & \text{if } x > 1. \end{cases}$$

*Then* $\Psi \circ \psi^{-1}$ *is regularly varying at infinity with index*

$$\gamma = \begin{cases} 0 & \text{if } \beta - \alpha \leq -1, \\ (\beta - \alpha + 1)/\alpha & \text{if } \beta - \alpha > -1. \end{cases}$$

*Furthermore, let* $X$ *be a random variable and let* $a > 0$. *Then*

$$\int_a^\infty \left(\sum_{n \geq 1} \phi(n) P(|X| > x\psi(n))\right) dx < \infty \text{ if and only if } E\left(|X|\Psi(\psi^{-1}(|X|))\right) < \infty.$$

**Lemma 3.3** (Shao (2000)). *Let* $\{X_i; i \geq 1\}$ *be a sequence of NA random variables with* $EX_i = 0, i \geq 1$.

(a) *If* $E|X_i|^p < \infty$ *for some* $p \geq 2$ *and all* $i \geq 1$, *then there exists a constant* $A_p > 0$ *such that for all* $n \geq 1$

$$E\left(\max_{1\leq k\leq n} |\sum_{i=1}^{k} X_i|^p\right) \leq A_p \left\{ \left(\sum_{i=1}^{n} EX_i^2\right)^{\frac{p}{2}} + \sum_{i=1}^{n} E|X_i|^p \right\}.$$

(b) *Set* $B_n^2 = \sum_{i=1}^{n} EX_i^2, n \geq 1$. *If the* $\{X_i\}$ *have finite variance, then for all* $\xi > 0, \eta > 0$, *and* $0 < \lambda < 1$, *we have for all* $n \geq 1$

$$P\left(\max_{1\leq k\leq n} |\sum_{i=1}^{k} X_i| \geq \xi\right) \leq 2P\left(\max_{1\leq k\leq n} |X_k| > \eta\right) + \frac{2}{1-\lambda} \exp\left\{-\frac{\xi^2 \lambda}{2(\eta\xi + B_n^2)}\right\}.$$

**Lemma 3.4.** *Let* $\alpha \geq -1, \beta \geq 0$, *and let* $\{Y_n; n \geq 1\}$ *be a nondecreasing sequence of nonnegative random variables. If*

$$\sum_{n\geq 1} n^\alpha P\left(Y_n > xn^\beta\right) < \infty \quad \text{for all} \quad x > 0, \tag{3.1}$$

*then*

$$\lim_{n\to\infty} P\left(Y_n > xn^\beta\right) = 0 \quad \text{for all} \quad x > 0. \tag{3.2}$$

*Proof.* The result is obvious if $\alpha \geq 0$ so assume that $-1 \leq \alpha < 0$. Note that $0 \leq Y_1 \leq Y_2 \leq \cdots \leq Y_n \leq \cdots$ and for all $x > 0$, recalling (3.1) we have

$$\begin{aligned}
\infty &> \sum_{n\geq 1} n^\alpha P\left(Y_n > xn^\beta\right) \\
&= \sum_{j\geq 0} \sum_{2^j \leq n < 2^{j+1}} n^\alpha P\left(Y_n > xn^\beta\right) \\
&\geq \sum_{j\geq 0} \sum_{2^j \leq n < 2^{j+1}} n^\alpha P\left(Y_{2^j} > x 2^{(j+1)\beta}\right) \\
&= \sum_{j\geq 0} P\left(Y_{2^j} > 2^\beta x 2^{j\beta}\right) \sum_{2^j \leq n < 2^{j+1}} n^\alpha \\
&\geq \sum_{j\geq 0} P\left(Y_{2^j} > 2^\beta x 2^{j\beta}\right) 2^j 2^{(j+1)\alpha} \\
&= 2^\alpha \sum_{j\geq 0} 2^{(\alpha+1)j} P\left(Y_{2^j} > 2^\beta x 2^{j\beta}\right)
\end{aligned}$$

which ensures since $\alpha + 1 \geq 0$ that

$$\lim_{j\to\infty} P\left(Y_{2^j} > x_1 2^{j\beta}\right) = 0 \quad \text{for all} \quad x_1 > 0. \tag{3.3}$$

Now for each $n \geq 1$, let $j_n \geq 0$ be such that $2^{j_n} \leq n < 2^{j_n+1}$; i.e., $j_n = [\log_2 n], n \geq 1$ where $\log_2$ denotes the logarithm to the base 2 and $[x]$ denotes the greatest integer in $x$ for $x \geq 0$. Then for all $x > 0$,

$$\begin{aligned} P\left(Y_n > xn^\beta\right) &\leq P\left(Y_{2^{j_n+1}} > x2^{j_n\beta}\right) \\ &= P\left(Y_{2^{j_n+1}} > x2^{-\beta}2^{(j_n+1)\beta}\right) \\ &\to 0 \quad \text{as } n \to \infty \end{aligned}$$

by (3.3) thereby proving (3.2). $\square$

*Proof of Theorem 2.1.* The equivalence between (i) and (ii) follows by arguing as in Remark 1.1. We now prove the implication ((ii)⇒(iii)). Note that (ii) implies

$$\sum_{n \geq 1} n^{r-2} P\left(\max_{1 \leq k \leq n} |S_k| > x^q n^{\frac{1}{p}}\right) < \infty \quad \text{for all } x > 0. \tag{3.4}$$

Let $\alpha = r - 2$. Then $\alpha > -1$ since $r > 1$ and (3.1) ensures via Lemma 3.4 with $\beta = 1/p$ and $Y_n := \max_{1 \leq k \leq n} |S_k|$, $n \geq 1$ that $P\left(\max_{1 \leq k \leq n} |S_k| > x^q n^{1/p}\right) \to 0$ for all $x > 0$ as $n \to \infty$, and it follows from Lemma 3.1 that for all sufficiently large $n$ and all $x > 0$

$$nP\left(|X| > x^q n^{\frac{1}{p}}\right) \leq CP\left(\max_{1 \leq k \leq n} |S_k| > 2^{-1} x^q n^{\frac{1}{p}}\right). \tag{3.5}$$

Now (ii) and (3.5) yield for all $\epsilon > 0$

$$\int_\epsilon^\infty \sum_{n \geq 1} n^{r-1} P\left(|X|^{\frac{1}{q}} > xn^{\frac{1}{pq}}\right) dx < \infty,$$

which, on account of Lemma 3.2, is equivalent to

$$E\left[|X|^{\frac{1}{q}} \Psi(\psi^{-1}(|X|^{\frac{1}{q}}))\right] < \infty. \tag{3.6}$$

Here $\psi(x) = x^{\frac{1}{pq}}$, $\phi(x) = x^{r-1}$, $x \geq 0$, and (3.6) is equivalent to

$$\begin{cases} E|X|^{\frac{1}{q}} < \infty & \text{if } q < \frac{1}{pr}, \\ E[|X|^{pr} \log^+ |X|] < \infty & \text{if } q = \frac{1}{pr}, \\ E|X|^{pr} < \infty & \text{if } q > \frac{1}{pr}. \end{cases} \tag{3.7}$$

Next, we prove $EX = 0$ when $1 \leq p < 2$.

By the implication ((iii)⇒(ii)) (which is established below), (3.4) implies that for $1 \leq p < 2$

$$\int_\epsilon^\infty \sum_{n \geq 1} n^{r-2} P\left(\max_{1 \leq k \leq n} |S_k - ES_k| > x^q n^{\frac{1}{p}}\right) dx < \infty \quad \text{for all } \epsilon > 0,$$

which, together with (ii) and $r > 1$, ensures that $EX = 0$ when $1 \leq p < 2$.

Now, we prove the implication ((iii)⇒(ii)). Note that

$$\int_\epsilon^\infty f(x^q)\, dx \leq \sum_{n \geq 1} n^{r-2} P\left(\max_{1 \leq k \leq n} |S_k| > \epsilon^q n^{\frac{1}{p}}\right)$$
$$+ \sum_{n \geq 1} n^{r-2} \int_1^\infty P\left(\max_{1 \leq k \leq n} |S_k| > x^q n^{\frac{1}{p}}\right) dx.$$

Thus, it suffices to show that

$$\sum_{n \geq 1} n^{r-2} P\left(\max_{1 \leq k \leq n} |S_k| > \epsilon^q n^{\frac{1}{p}}\right) < \infty \tag{3.8}$$

and

$$\sum_{n \geq 1} n^{r-2} \int_1^\infty P\left(\max_{1 \leq k \leq n} |S_k| > x^q n^{1/p}\right) dx < \infty. \tag{3.9}$$

We prove only (3.9); the proof of (3.8) is analogous.

Choose $0 < \alpha < q$, $1/pr < \beta < 1/p$. Set for all $n \geq 1$ and $1 \leq i \leq n$

$$X_{ni}(1) = -x^\alpha n^\beta I\left(X_i < -x^\alpha n^\beta\right) + X_i I\left(|X_i| \leq x^\alpha n^\beta\right) + x^\alpha n^\beta I\left(X_i > x^\alpha n^\beta\right),$$
$$X_{ni}(2) = \left(X_i - x^\alpha n^\beta\right) I\left(x^\alpha n^\beta < X_i < x^q n^{1/p}/(4N)\right),$$
$$X_{ni}(3) = \left(X_i + x^\alpha n^\beta\right) I\left(-x^\alpha n^\beta > X_i > -x^q n^{\frac{1}{p}}/(4N)\right),$$
$$X_{ni}(4) = \left(X_i - x^\alpha n^\beta\right) I\left(X_i \geq x^q n^{\frac{1}{p}}/(4N)\right) + \left(X_i + x^\alpha n^\beta\right) I(X_i \leq -x^q n^{\frac{1}{p}}/(4N)),$$

where $N$ is a large positive integer which will be specified later. Then

$$\sum_{n \geq 1} n^{r-2} \int_1^\infty P\left(\max_{1 \leq k \leq n} |S_k| > x^q n^{\frac{1}{p}}\right) dx$$
$$\leq \sum_{n \geq 1} n^{r-2} \int_1^\infty P\left(\max_{1 \leq k \leq n} \left|\sum_{i=1}^k X_{ni}(1)\right| > \frac{x^q n^{1/p}}{4}\right) dx$$
$$+ \sum_{n \geq 1} n^{r-2} \int_1^\infty P\left(\max_{1 \leq k \leq n} \left|\sum_{i=1}^k X_{ni}(2)\right| > \frac{x^q n^{1/p}}{4}\right) dx$$
$$+ \sum_{n \geq 1} n^{r-2} \int_1^\infty P\left(\max_{1 \leq k \leq n} \left|\sum_{i=1}^k X_{ni}(3)\right| > \frac{x^q n^{1/p}}{4}\right) dx$$
$$+ \sum_{n \geq 1} n^{r-2} \int_1^\infty P\left(\max_{1 \leq k \leq n} \left|\sum_{i=1}^k X_{ni}(4)\right| > \frac{x^q n^{1/p}}{4}\right) dx$$
$$=: I_1 + I_2 + I_3 + I_4.$$

By applying Lemma 3.2, (iii) yields

$$I_4 \leq \sum_{n \geq 1} n^{r-1} \int_1^\infty P\left(|X| \geq \frac{x^q n^{1/p}}{4N}\right) dx < \infty.$$

From the definition of $X_{ni}(2)$, we know that $X_{ni}(2) > 0$. Taking $N > \max\left\{\frac{r-1}{\beta pr - 1}, \frac{1}{\alpha pr}\right\}$, by the property of NA, we have

$$\begin{aligned}
I_2 &= \sum_{n \geq 1} n^{r-2} \int_1^\infty P\left(\max_{1 \leq k \leq n} \left|\sum_{i=1}^k X_{ni}(2)\right| > \frac{x^q n^{1/p}}{4}\right) dx \\
&\leq \sum_{n \geq 1} n^{r-2} \int_1^\infty P\left(\begin{array}{c} \text{there are at least } N \text{ indices } i \in [1, n] \\ \text{such that } X_{ni}(2) \neq 0 \end{array}\right) dx \\
&\leq \sum_{n \geq 1} n^{r-2} \int_1^\infty \left(\sum_{i=1}^n P(X_i > x^\alpha n^\beta)\right)^N dx \\
&\leq \sum_{n=1}^\infty n^{r-2-(\beta pr - 1)N} (E|X|^{pr})^N \int_1^\infty x^{-\alpha pr N} dx < \infty.
\end{aligned}$$

Similarly, $X_{ni}(3) < 0$ and $I_3 < \infty$.

In order to estimate $I_1$, we first verify that

$$\frac{\max_{1 \leq k \leq n} \left|\sum_{i=1}^k EX_{ni}(1)\right|}{x^q n^{1/p}} \to 0. \tag{3.10}$$

Note that (iii) implies $E|X|^{pr} < \infty$ and $E|X|^p < \infty$. When $1 \leq p < 2$, by $EX = 0$, and noticing that $1 - \beta(pr - 1) - (1/p) < 0$ we have

$$\begin{aligned}
&\frac{\max_{1 \leq k \leq n} |\sum_{i=1}^k EX_{ni}(1)|}{x^q n^{1/p}} \\
&\leq \frac{2nE|X|I(|X| > x^\alpha n^\beta)}{x^q n^{1/p}} \\
&\leq 2x^{-\alpha(pr-1)-q} n^{1-\beta(pr-1)-\frac{1}{p}} E|X|^{pr} I\left(|X| > x^\alpha n^\beta\right) \\
&\to 0 \text{ as } n \to \infty.
\end{aligned} \tag{3.11}$$

When $0 < p < 1$, if $pr > 1$, we have

$$\begin{aligned}
&\frac{\max_{1 \leq k \leq n} |\sum_{i=1}^k EX_{ni}(1)|}{x^q n^{1/p}} \\
&\leq \frac{nx^\alpha n^\beta P(|X| > x^\alpha n^\beta)}{x^q n^{1/p}} + \frac{nE|X|I(|X| \leq x^\alpha n^\beta)}{x^q n^{1/p}} \\
&\leq x^{-\alpha(pr-1)-q} n^{1-\beta(pr-1)-(1/p)} E|X|^{pr} I(|X| > x^\alpha n^\beta) \\
&\quad + x^{-q} n^{-(1/p-1)} E|X| \\
&\to 0 \text{ as } n \to \infty,
\end{aligned} \tag{3.12}$$

and if $0 < pr < 1$, we have

$$\frac{\max_{1 \leq k \leq n} |\sum_{i=1}^{k} EX_{ni}(1)|}{x^q n^{1/p}} \leq 2x^{-\alpha(pr-1)-q} n^{1-\beta(pr-1)-(1/p)} E|X|^{pr} \quad (3.13)$$
$$\to 0 \quad \text{as} \quad n \to \infty.$$

Hence (3.10) follows from (3.11)-(3.13). Thus, to prove $I_1 < \infty$, we need only show that

$$I_1^* =: \sum_{n \geq 1} n^{r-2} \int_1^\infty P\left(\max_{1 \leq k \leq n} \left|\sum_{i=1}^k (X_{ni}(1) - EX_{ni}(1))\right| > \frac{x^q n^{1/p}}{5}\right) dx < \infty.$$

Note that $\{X_{ni}(1) - EX_{ni}(1); 1 \leq i \leq n, n \geq 1\}$ is an array of rowwise NA random variables from the definition of $X_{ni}(1)$. Here we will apply Lemma 3.3(a). Choosing

$$M > \max\left\{2, \ pr, \ \frac{1}{q}, \ \frac{2p(r-1)}{2-p}, \ \frac{1}{(q-\alpha) + \frac{\alpha pr}{2}}, \ \frac{r-1}{(\frac{1}{p} - \beta) + \frac{\beta pr - 1}{2}}, \ \frac{1 - \alpha pr}{q - \alpha}\right\}, \quad (3.14)$$

we have

$$I_1^* \leq C \sum_{n \geq 1} n^{r-2-\frac{M}{p}} \int_1^\infty x^{-Mq} \left\{\left(\sum_{i=1}^n EX_{ni}^2(1)\right)^{M/2} + \sum_{i=1}^n E|X_{ni}(1)|^M\right\} dx$$
$$=: I_{11} + I_{12}.$$

We observe that

$$I_{11} \leq C \sum_{n \geq 1} n^{r-2-\frac{M}{p}+\beta M + \frac{M}{2}} \int_1^\infty x^{-M(q-\alpha)} \left(P(|X| > x^\alpha n^\beta)\right)^{M/2} dx$$
$$+ C \sum_{n \geq 1} n^{r-2-\frac{M}{p}+\frac{M}{2}} \int_1^\infty x^{-Mq} \left(EX^2 I(|X| \leq x^\alpha n^\beta)\right)^{M/2} dx$$
$$=: I_{11}^{(1)} + I_{11}^{(2)}.$$

From (3.14) and $E|X|^{pr} < \infty$ we have

$$I_{11}^{(1)} \leq C \sum_{n \geq 1} n^{r-2-[(\frac{1}{p}-\beta)+\frac{\beta pr-1}{2}]M} (E|X|^{pr})^{M/2} \int_1^\infty x^{-[(q-\alpha)+\frac{\alpha pr}{2}]M} dx < \infty.$$

As for $I_{11}^{(2)}$, if $pr \geq 2$, we have

$$I_{11}^{(2)} \leq C \sum_{n \geq 1} n^{r-2-(\frac{1}{p}-\frac{1}{2})M} (EX^2)^{M/2} \int_1^\infty x^{-qM} dx < \infty,$$

and if $0 < pr < 2$, we have

$$I_{11}^{(2)} \leq C \sum_{n \geq 1} n^{r-2-[(\frac{1}{p}-\beta)+\frac{\beta pr-1}{2}]M}(E|X|^{pr})^{M/2} \int_1^\infty x^{-[(q-\alpha)+\frac{\alpha pr}{2}]M} dx < \infty.$$

Therefore, $I_{11} < \infty$.

As for $I_{12}$, from (3.14) and $E|X|^{pr} < \infty$ we obtain that

$$\begin{aligned}
I_{12} &\leq C \sum_{n \geq 1} n^{r-1-(\frac{1}{p}-\beta)M} \int_1^\infty x^{-M(q-\alpha)} P(|X| > x^\alpha n^\beta) dx \\
&\quad + C \sum_{n \geq 1} n^{r-1-\frac{M}{p}} \int_1^\infty x^{-Mq} E|X|^M I(|X| \leq x^\alpha n^\beta) dx \\
&\leq C \sum_{n \geq 1} n^{r-1-(\frac{1}{p}-\beta)M-\beta pr}(E|X|^{pr}) \int_1^\infty x^{-M(q-\alpha)-\alpha pr} dx \\
&< \infty.
\end{aligned}$$

This establishes that $I_1^* < \infty$ thereby completing the proof of Theorem 2.1. □

*Proof of Theorem 2.2.* The equivalence between (i) and (ii) follows by again arguing as in Remark 1.1, *mutatis mutandis*. We only give the proof of the implication ((iii)⇒(ii)) and leave the proof of the implication ((ii)⇒(iii)), which is similar to the argument for ((ii)⇒(iii)) in Theorem 2.1, to the reader.

Choose $0 < \alpha < q$. Set for all $n \geq 1$ and $1 \leq i \leq n$

$$\begin{aligned}
Y_{ni}(1) &= -x^\alpha(n/\log n)^{1/2} I\left(X_i < -x^\alpha(n/\log n)^{1/2}\right) \\
&\quad + X_i I\left(|X_i| \leq x^\alpha(n/\log n)^{1/2}\right) \\
&\quad + x^\alpha(n/\log n)^{1/2} I\left(X_i > x^\alpha(n/\log n)^{1/2}\right), \\
Y_{ni}(2) &= \left(X_i - x^\alpha(n/\log n)^{1/2}\right) \\
&\quad I\left(x^\alpha(n/\log n)^{1/2} < X_i < x^q(n\log n)^{1/2}/(4N)\right), \\
Y_{ni}(3) &= \left(X_i + x^\alpha(n/\log n)^{1/2}\right) \\
&\quad I\left(-x^\alpha(n/\log n)^{1/2} > X_i > -x^q(n\log n)^{1/2}/(4N)\right), \\
Y_{ni}(4) &= \left(X_i - x^\alpha(n/\log n)^{1/2}\right) I\left(X_i \geq x^q(n\log n)^{1/2}/(4N)\right) \\
&\quad + \left(X_i + x^\alpha(n/\log n)^{1/2}\right) I\left(X_i \leq -x^q(n\log n)^{1/2}/(4N)\right),
\end{aligned}$$

where $N$ is a large positive integer which will be specified later. Then

$$\sum_{n\geq 3} n^{r-2} \int_\epsilon^\infty P\left(\max_{1\leq k\leq n} |S_k| > x^q(n\log n)^{1/2}\right) dx$$

$$\leq \sum_{n\geq 3} n^{r-2} \int_\epsilon^\infty P\left(\max_{1\leq k\leq n} \left|\sum_{i=1}^k Y_{ni}(1)\right| > \frac{x^q(n\log n)^{1/2}}{4}\right) dx$$

$$+ \sum_{n\geq 3} n^{r-2} \int_\epsilon^\infty P\left(\max_{1\leq k\leq n} \left|\sum_{i=1}^k Y_{ni}(2)\right| > \frac{x^q(n\log n)^{1/2}}{4}\right) dx$$

$$+ \sum_{n\geq 3} n^{r-2} \int_\epsilon^\infty P\left(\max_{1\leq k\leq n} \left|\sum_{i=1}^k Y_{ni}(3)\right| > \frac{x^q(n\log n)^{1/2}}{4}\right) dx$$

$$+ \sum_{n\geq 3} n^{r-2} \int_\epsilon^\infty P\left(\max_{1\leq k\leq n} \left|\sum_{i=1}^k Y_{ni}(4)\right| > \frac{x^q(n\log n)^{1/2}}{4}\right) dx$$

$$=: J_1 + J_2 + J_3 + J_4.$$

By applying Lemma 3.2, as in the proof of ((iii)⇒(ii)) in Theorem 2.1, (iii) yields

$$J_4 \leq \sum_{n\geq 3} n^{r-1} \int_\epsilon^\infty P\left(|X| \geq \frac{x^q(n\log n)^{1/2}}{4N}\right) dx < \infty.$$

Following the line of argument as for bounding $I_2$ in the proof of Theorem 2.1, we have

$$J_2 \leq \sum_{n\geq 3} n^{r-2} \int_\epsilon^\infty \left(\sum_{i=1}^n P(X_i > x^\alpha(n/\log n)^{1/2})\right)^N dx < \infty.$$

Similarly, $Y_{ni}(3) < 0$ and $J_3 < \infty$.

Note that (ii) implies $EX^2 < \infty$. Therefore, by $EX = 0$ we have

$$\frac{\max_{1\leq k\leq n} |\sum_{i=1}^k EY_{ni}(1)|}{x^q(n\log n)^{1/2}} \leq \frac{2nE|X|I(|X| > x^\alpha(n/\log n)^{1/2})}{x^q(n\log n)^{1/2}} \to 0 \text{ as } n \to \infty.$$

Thus, to prove $J_1 < \infty$, it suffices to show that

$$J_1^* =: \sum_{n\geq 3} n^{r-2} \int_\epsilon^\infty P\left(\max_{1\leq k\leq n} \left|\sum_{i=1}^k (Y_{ni}(1) - EY_{ni}(1))\right| > \frac{x^q(n\log n)^{1/2}}{5}\right) dx < \infty.$$

Here we will apply Lemma 3.3(b). Taking $\xi = x^q(n\log n)^{1/2}/5$, $\eta = 2x^\alpha(n/\log n)^{1/2}$, $\lambda = 1/2$. Note that $\max_{1\leq k\leq n} |Y_{ni}(1) - EY_{ni}(1)| \leq \eta$, $B_n^2 = \sum_{i=1}^n (Y_{ni}(1) - EY_{ni}(1))^2 \leq 2nEX^2$. Hence, according to Lemma 3.3(b), for $x \geq (20EX^2)^{1/(\alpha+q)}$ we have

$$P\left(\max_{1\leq k\leq n} \left|\sum_{i=1}^k (Y_{ni}(1) - EY_{ni}(1))\right| > \frac{x^q(n\log n)^{1/2}}{5}\right) \leq 4\exp\left\{-\frac{x^{q-\alpha}\log n}{2}\right\}.$$

Now, choosing $\epsilon > \max\{[2(r-1)]^{1/(q-\alpha)},\ (20EX^2)^{1/(\alpha+q)}\}$ we get

$$J_1^* \leq 4 \sum_{n \geq 3} n^{r-2} \int_\epsilon^\infty \exp\left\{-\frac{x^{q-\alpha}\log n}{2}\right\} dx < \infty.$$

This establishes that $J_1 < \infty$ thereby completing the proof of Theorem 2.2. □